\documentclass{elsarticle}
\usepackage{amsfonts}

\usepackage{graphicx}
\usepackage{pstricks}
\usepackage{amsmath}
\usepackage{amsxtra}
\newtheorem{theorem}{Theorem}[section]
\newtheorem{lemma}[theorem]{Lemma}
\newproof{proof}[theorem]{Proof}
\newtheorem{proposition}[theorem]{Proposition}
\newtheorem{corollary}[theorem]{Corollary}

\newtheorem{problem}[theorem]{Problem}

\newtheorem{example}[theorem]{Example}
\newtheorem{algorithm}[theorem]{Algorithm}

\numberwithin{equation}{section}

\renewcommand{\phi}{\varphi}

\newcommand{\V}{\mathcal{V}}

\renewcommand{\emph}{\textsl}

\begin{document}

\begin{frontmatter}

\title{Non-extremal Sextic Moment Problems}

\author{Ra\'{u}l E. Curto\footnote{The first named author was partially supported by NSF Grants
DMS-0801168 and DMS-1302666.}}
\address{Department of Mathematics, The University of Iowa, Iowa City, Iowa
52242}
\ead{raul-curto@uiowa.edu}
\ead[url]{http://www.math.uiowa.edu/\symbol{126}rcurto/}

\author{Seonguk Yoo\footnote{The second named author was supported by the Brain Korea 21 Program of the National Research Foundation of Korea.}}
\address{Department of Mathematics, Inha University, Incheon 402-751, Korea}
\ead{seyoo73@gmail.com}

\begin{abstract}
Positive semidefiniteness, recursiveness,  and the variety condition of a moment
matrix are necessary and sufficient conditions to solve the
quadratic and quartic moment problems.\ Also, positive semidefiniteness, combined with another
necessary condition, consistency, is a sufficient condition in the case of
\textit{extremal} moment problems, i.e., when the rank of the moment matrix (denoted by $r$) and the cardinality of the
associated algebraic variety (denoted by $v$) are equal. \ However,  these conditions
are not sufficient for non-extremal sextic or higher-order truncated moment
problems. \

In this paper we settle three key instances of the \textit{non}-extremal (i.e., $r<v$) sextic moment problem, as follows: when $r=7$, positive semidefiniteness, consistency and the variety condition guarantee the existence of a $7$-atomic representing measure; when $r=8$ we construct two determining algorithms, corresponding to the cases $v=9$ and $v=+ \infty$. \ To accomplish this, we generalize the rank-reduction technique developed in previous work, where we solved the nonsingular quartic moment problem and found an explicit way to build a representing measure.
\end{abstract}


\end{frontmatter}


\section{Introduction}\label{intro}

Given a collection of real numbers $\beta \equiv\beta^{(2n)}= \{  \beta_{00},$ $\beta_{10},$  $\beta_{01},\cdots, \beta_{2n,0},\beta_{2n-1,1},\newline \cdots, \beta_{1,2n-1},\beta_{0,2n} \}$, the \textit{truncated real moment problem} (TMP) consists of finding a positive Borel measure $\mu$ supported in the real plane $\mathbb{R}^2$ such that
\begin{eqnarray*}
\beta_{ij}=\int x^i y^j \,\, d\mu \,\,\,(i,j\in \mathbb Z_+, \ 0\leq i+j\leq 2n).
\end{eqnarray*}
The collection $\beta$ is called a \textit{truncated moment sequence} (of order $2n$) and $\mu $ is called a 
\textit{representing measure} for $\beta $.

Naturally associated with each TMP is a moment matrix $%
\mathcal{M}(n)\equiv \mathcal{M}(n)(\beta )$, defined by 
$$
\mathcal{M}(n)(\beta):=(\beta_{\textbf{ i} +\textbf{j} })_{\textbf{i}, \; \textbf{j} \;\in \mathbb Z^2_+: |\textbf{i}|, |\textbf{j}|\leq n}.
$$
In order to define a functional calculus for the columns of $\mathcal{M}(n)$, we label the columns of  $\mathcal M{(n)}$ with the following lexicographical order: ${\it 1},X,Y,X^2,X Y,Y^2, \cdots,X^n,\cdots,Y^n$. \ Let $M[i,j]$ be the Hankel matrix of size  $(i+1)\times(j+1)$ whose entries are moments of order $i+j$, as follows:
$$
M[i,j]:=\left(
               \begin{array}{cccc}
                  \beta_{i+j,0} & \beta_{i-1+j,1}&\cdots & \beta_{i,j} \\
                 \beta_{i-1+j,1} & \beta_{i-2+j,2}&\cdots & \beta_{i-1,1+j} \\
\vdots & \ddots&\ddots&\vdots\\
                  \beta_{j,i} & \beta_{j-1, i+1}  &\cdots & \beta_{0,i+j}
                   \end{array}
              \right).
$$

Note that $\mathcal{M}(n) = (M(i,j))_{i,j=1,\cdots,n}$, and that 
$$
\mathcal{M}(n+1)=\left(
\begin{array}{cc}
\mathcal{M}(n) & B \\
B^{\ast } & M(n+1,n+1)%
\end{array}
\right),
$$ 
where $B$ is the block matrix $(M(i,n+1))_{i=1,\cdots,n}$. \ The matrix $\mathcal{M}(n)$ detects the positivity of the \textit{Riesz functional} $\Lambda: p \mapsto \sum_{ij}a_{ij}\beta _{ij}\;\;(p(x,y)\equiv
\sum_{ij}a_{ij} x^{i}y^{j})$ on the cone generated by the collection $%
\{p^2:p\in \mathbb{R}[x,y]\}$. \ 

In addition to its importance
for applications, a complete solution of TMP would readily lead to a
solution of the \textit{full} moment problem, via a weak-* convergence argument, as
shown by J. Stochel \cite{Sto2}. \ While we primarily focus on truncated moment problems, the full moment problem (in one or several variables) has been widely studied; see, for example, \cite{BaTe}, \cite{blekherman}, \cite{Dem}, \cite{KuMa}, \cite{Lau2}, \cite{Lau3}, \cite{PoSc}, \cite{Pu1}, \cite{PuSch}, \cite{PuSchm}, \cite{PuVa1}, \cite{PuVa2}, \cite{Rie}, \cite{Sche1}, \cite{Sche2}, \cite{Sch1}, \cite{Sch2}, \cite{Sch4}, \cite{StSz1}, \cite{polar}, \cite{Vas2}.

Building on previous work for the case of \textit{real} moments, several years ago the first named author and L. Fialkow introduced in 
\cite{tcmp1}, \cite{tcmp2} and \cite{tcmp3} an approach to TMP based on
matrix positivity and extension, combined with a new ``functional calculus''
for the columns of $\mathcal{M}(n)$. \ This allowed them to show that TMP is soluble in the following cases:

\noindent (i) TMP is of \textit{flat data} type \cite{tcmp1}, i.e., $\operatorname{rank}\;\mathcal{M}(n)=\operatorname{rank}\;\mathcal{M}(n-1)$;

\noindent(ii) the columns ${\it 1},X,Y$ are linearly dependent \cite[Theorem 2.1]{tcmp2};

\noindent(iii) $\mathcal{M}(n)$ is singular and subordinate to conics \cite{tcmp5}, \cite%
{tcmp6}, \cite{tcmp7}, \cite{tcmp9};

\noindent(iv) $\mathcal{M}(n)$ admits a rank-preserving moment matrix extension $\mathcal{M}(n+1)$, i.e., an extension $\mathcal{M}(n+1)$ which is flat 
\cite{tcmp10};

\noindent(v) $\mathcal{M}(n)$ is extremal, i.e., $\operatorname{rank}\;\mathcal{M}(n)=\operatorname{card} \;%
\mathcal{V}(\beta ^{(2n)})$, where $\mathcal{V}(\beta )\equiv \mathcal{V}%
(\beta ^{(2n)})$ is the algebraic variety of $\beta $ \cite{tcmp11}.

\noindent(vi) $\mathcal{M}(n)$ is \textit{recursively determinate}, that is,  if $\mathcal{M}(n)$ has only column dependence relations of the form
\begin{eqnarray*}
&&X^n = p(X,Y) ~~~ (p\in \mathcal P_{n-1});\\ 
&&Y^m = q(X,Y) ~~~(q\in \mathcal P_m, ~ q \textrm{ has no $y^m$ term, } m \leq n),
\end{eqnarray*}
where $\mathcal P_k$ denotes the subspace of polynomials in $\mathbb R[x,y]$ whose degree is less than or equal to $k$ \cite{tcmp13}.

The common feature of the above mentioned cases is
the presence, at the level of the column space $\mathcal{C}_{\mathcal{M}(n)}$, of algebraic conditions
implied by the existence of a representing measure with support in a proper
real algebraic subset of the plane. \ However, the chief attraction
of the truncated moment problem (TMP) is its naturalness: since the data set is
finite, we can apply ``finite'' techniques, grounded in finite dimensional
operator theory, linear algebra, and algebraic geometry, to develop
algorithms for explicitly computing finitely atomic representing measures. \ 

Consistent with this view, in this paper we solve the non-extremal sextic moment problem, as follows. \ Without loss of generality, one assumes that $\mathcal{M}(2)$ is invertible, and that $\mathcal{M}(3)$ is not a flat extension of $\mathcal{M}(2)$; that is, $r \ge 7$. \ When $r=7$, positive semidefiniteness, consistency and the variety condition guarantee the existence of a $7$-atomic representing measure; when $r=8$ we construct two determining algorithms, corresponding to the cases $v=9$ and $v=+ \infty$. \ To accomplish this, we generalize the rank-reduction technique developed in \cite{CuYoo3}, where we solved the nonsingular quartic moment problem and found an explicit way to build a representing measure.  


\subsection{\textbf{Necessary Conditions}\label{NCS}}
In order to introduce basic necessary conditions for the existence of a measure, let $\mu$ be a representing measure of $\beta$. First, we can compute  that
\begin{eqnarray*}
0 \leq \int| p(x,y)|^2 \, d\mu
 = \sum_{i,j,k,\ell} a_{ij}a_{k\ell} \int x^{i+k} y^{j+\ell} \, d\mu
 = \sum_{i,j,k,\ell} a_{ij} a_{k\ell} \beta_{i+k} \beta_{j+\ell}
\end{eqnarray*}
which is equivalent to the condition $\mathcal{M}(n) \geq 0$. \
We next define an assignment from $\mathcal P_n$ to $\mathcal C_{\mathcal{M}(n)}$; given a polynomial $p(x,y) \equiv \sum_{ij}a_{ij}x^{i}y^{j}$, we let $p(X,Y):=\sum_{ij}a_{ij}X^{i}Y^{j}$ (so that $p(X,Y
)\in \mathcal{C}_{\mathcal{M}(n)}$), which is the so-called ``functional calculus." \
We also let $\mathcal{Z}(p)$ denote the zero set of $p$ and define the \textit{algebraic variety} of $\beta $ by
\begin{equation} \label{variety}
\mathcal{V} \equiv \mathcal{V}(\beta ):=\bigcap {}_{p\, (X,Y)=0,\, \deg \, p\, \leq n}\; \mathcal{Z}%
(p).\
\end{equation}%
If $\widehat{p}$ denotes the column vector of coefficients of $p$, then we see at once that $p(X,Y)=\mathcal{M}(n)\widehat{p}$, that is, $p(X,Y)=0$ if and only if $\widehat{p} \in \operatorname{ker}\; \mathcal{M}(n)$. \ It follows that $\operatorname{supp}\; \mu \subseteq \mathcal{V}(\beta )$ and
$r:=\operatorname{rank}\; \mathcal{M}(n)\leq \operatorname{card}\;\operatorname{supp}\;\mu \leq v:= \operatorname{card}\; \V$ \cite{tcmp3}. \ Thus, $r \leq v$ is a necessary condition for solubility, referred to as the \textit{variety condition}. \
In addition, if $p$ is any polynomial of degree at most $2n$ such that $p|_{%
\mathcal{V}}\equiv 0$, then the Riesz functional $\Lambda $ must satisfy $%
\Lambda (p):=\int p \, d\mu =0$; this is the so-called \textit{consistency} of the moment sequence, and is also a necessary condition for solubility \cite{tcmp11}.\

Positive semidefiniteness, recursiveness and the variety condition solve the quadratic and quartic moment problems (see \cite{tcmp1}, \cite{tcmp6}, \cite{FiaNie}). \ Moreover, the main result in \cite{tcmp11} establishes that the preceding three conditions together with consistency are sufficient in the case of {\it extremal} moment problems (i.e., $r=v$). \ In \cite{tcmp11}, the authors also showed that consistency cannot be replaced by the weaker condition that $\mathcal{M}(n)$ is \textit{recursively generated}, (RG), that is, if  $p(X,Y)=0$, then $(p\, q)(X,Y)=0$, for each polynomial $q$ with $\deg (p\, q)\leq n$.

Any \textit{singular} moment matrix $\mathcal{M}(n)$ must have at least one
linear column dependence relation. \ A result independently proven by \ H.M. M\"oller \cite{Moe} and C. Scheiderer
says that  the polynomials associated to the column relations  generate a real radical ideal $\mathcal{I}$  whenever $\mathcal{M}(n)\geq 0$ (cf. \cite[Subsection 5.1, p.203]{Lau3} and \cite{Lau2}). \
Through this fact, we can exploit some results from algebraic geometry for the study of TMP; for example, we have applied the Division Algorithm for multivariable polynomials to obtain a structure theorem that plays an important role in getting the main results in \cite{CuYoo2}. \


\subsection{\textbf{Flat Extensions}}

We recall that if $\mathcal{M}(n)$ admits a flat extension $\mathcal M(n+1)$, then  $\beta$ has an $\operatorname{rank}\; \mathcal{M}(n)$-atomic measure. \
This result is called the Flat Extension Theorem and is the most general solution to TMP, even though the actual construction of an extension is sometimes not feasible for high-order TMP. \
This theorem will be used implicitly in the proof of the first main result, Theorem \ref{main1}. \ For the reader's convenience, we recall that the general form of a flat extension of a positive semidefinite matrix $A$ is given by

$$
\left(
\begin{array}{cc}
A & AW \\
W^{\ast }A & W^{\ast }AW%
\end{array}
\right)
$$
for some matrix $W$. \ Note that while flat extensions have a simple algebraic structure, generating flat extensions of a moment matrix $A \equiv \mathcal{M}(n)$ requires verification that the blocks $AW$ and $W^{\ast }AW$ satisfy the relevant Hankel properties.\
In other words, an extension $\mathcal M(n+1)$ must be positive semidefinite while maintaining the moment matrix structure.


\subsection{\textbf{Centrality of  Extremal Moment Problems}}

The results in \cite{tcmp1} and \cite{tcmp3} show that any  \textit{soluble} TMP with a finite algebraic variety must have a moment matrix extension which is extremal (see also \cite{Fia08}).\ We thus need to find the minimal integer $k$ satisfying $\operatorname{rank}\;\mathcal M (n+k)=\operatorname{rank}\;\mathcal M (n+k+1).$\
However, the process is intricate; for example, even the latest version of  \textit{Mathematica} is unable to deliver the symbolic calculation needed to generate the most general extension $\mathcal M(3)$ of an invertible $\mathcal M(2)$.\
The existence criterion in \cite[Theorem 1.5]{tcmp3} provided an upper bound for the length of the extension sequence, which is  $k= 2n^2 + 6n + 6$; 
while the criterion sets a finite bound, it is definitely not sharp; for example, if $n=3$, then $k= 42$.\
However, for the cases with a finite variety, the following theorem gives  a significantly sharper bound.\
\begin{theorem}\cite{Fia08}\label{min-k}
 Suppose $v<\infty$. \ Then ${\beta}$ admits a representing measure if and only if
  $\mathcal M (n)({\beta})$ has a positive extension $\mathcal M (n+v-r+1)$ satisfying $\operatorname{rank}\;\mathcal M (n+v-r+1)\leq\operatorname{card}\;\mathcal V_{\mathcal M (n+v-r+1)}.$
\end{theorem}
The number $v-r$ is called the \textit{extremality gap}. \
This theorem says that for a moment sequence to have a representing measure, $\mathcal M (n)$ must have an extension sequence (with ranks possibly increasing) with maximal length $k= v-r$. \ For instance, if $\mathcal M (3)$ is singular with invertible $\mathcal M (2)$ and $v<\infty$, then  the maximal cardinality of a finite $\V$ is 9 and the minimal rank is 7 (if we assume that $\mathcal{M}(2)$ is invertible and that $\mathcal{M}(3)$ is not a flat extension of $\mathcal{M}(2)$). \ Thus, we get $k\leq  v-r \leq 9-7=2$ and it is sufficient  to  check extensions up to $\mathcal M (3+2+1)=\mathcal M (6)$.\


\subsection{\textbf{The Sextic Moment Problem}}
Using Theorem \ref{min-k}, and focusing on the values $r$, $v$ and the extremality gap $v-r$, we can classify all sextic truncated moment problems, as follows. \ Let us first denote $ r_n:=\operatorname{rank}\;\mathcal M (n)$ and $v_n:=\operatorname{card}\;\mathcal  V(\mathcal M (n))$.\
(Hereafter, we will use this notation throughout our presentation.)\
First, note that each extension must satisfy the variety condition, and so it is necessary for the following chain of inequalities to be true:
$$
{r_n\leq r_{n+1}\leq r_{n+2} \leq \cdots \leq v_{n+2} \leq v_{n+1}\leq v_n}.
$$
Since we have a complete solution to the quartic moment problem, we may always assume that the submatrix  $\mathcal M (2)$ in $\mathcal M(3)$ is invertible, that is, $\mathcal M(2)>0$. \
In Table $1$ we list all the possible sextic moment problems and provide some  information for each case; 
we notice that the case of $r_3=v_3=9$ cannot happen. \ For,  if the rank is 9, the moment matrix has  only one column relation, which means the algebraic variety is the graph of an algebraic curve in the plane, and therefore infinite.  \
\begin{table}[hptp]\label{table1}
\centering
\caption{\footnotesize{Classification of sextic moment problems according to $r$ and $v$}}
\smallskip
\begin{tabular}{|c|c|c|c|c|c|c|}
  \hline
  $r_3$ & $v_3$ & $v_3-r_3$ & Max Extension & & Solution Presented in\\
  \hline
  \hline
  7 &  7 &  0 & $\mathcal M (4)$ & extremal & \cite{CuYoo1}, \cite{CuYoo2} \\
 \hline
   7 & 8 & 1 & $\mathcal M (5)$ & & Theorem \ref{main1} below \\
 \hline

   7 & 9 & 2 & $\mathcal M (6)$ & & Theorem \ref{main1} below \\
  \hline
  7 & $\infty$ &  N/A & N/A & & Theorem \ref{main1} below \\
  \hline
  8 &  8 & 0 & $\mathcal M (4)$ & extremal& \cite{CuYoo2}\\
\hline
   8 & 9 & 1 & $\mathcal M (5)$ &  &  Algorithm \ref{alg15} \\
 \hline

 8 & $\infty$ &  N/A & N/A & & Algorithm \ref{Alg} \\
  \hline

9 & $\infty$ &  N/A & N/A & & \cite{Fia11} (particular cases) \\
  \hline

10 & $\infty$ &  N/A & N/A & & unknown \\
  \hline
\end{tabular}
\end{table}

In general, a TMP with infinite algebraic variety is much more difficult to study than one with a finite variety; we know that the atoms of a representing measure must lie on the graph of an algebraic curve, but we don't know exactly which points are in the support of the measure when we attempt to build a flat extension. \ This might be one of the reasons why, for  $n\geq 3$,  the nonsingular TMP and the TMP admitting a single column relation remain unsolved yet. \

Since a solution to  the extremal case of the sextic moment problem (i.e., when the extremality gap is zero) is provided in \cite{CuYoo2}, in this paper we focus on the cases with nonzero extremality gap.\

We conclude this subsection listing a key proposition which will be used in the proof of the main results.

\begin{proposition}\label{3-1-tcmp1} (cf. \cite[Proposition 3.1]{tcmp1}) \ Suppose $\mu$ is a representing measure for $\beta$. For $p\in \mathcal P_n$,
$$
\operatorname{supp}\; \mu \subseteq \mathcal Z(p)  \iff p(X,Y)=0.
$$
\end{proposition}


\subsection{\textbf{Statement of the Main Results}}
It is well-known that any rank-one positive semidefinite matrix must be of the form ${\bf x} {\bf x}^\ast$ for some nonzero vector ${\bf x}\in \mathbb C^n$.\
 Also, a  positive matrix $A$ can be written as a sum:
\begin{eqnarray*}
A= \sum_{i=1}^{k}{\bf x}_i {\bf x}_i^\ast,
\end{eqnarray*}
for some nonzero  vectors ${\bf x}_i \in \mathbb C^n$ for $i=1,\ldots,k$.\
The minimum number of summands is the rank of $A$.\ On the other
hand, if $\mathcal{M}(n)$ has an $r$-atomic measure $\mu \equiv \sum_{i=1}^{r} \rho_i \delta_{(x_i,y_i)}$, then we may write $\mathcal{M}(n)$
as
\begin{eqnarray*}
\mathcal{M}(n)= \sum_{i=1}^{r} \rho_i {\bf v}_i {\bf v}_i^\ast,
\end{eqnarray*}
where the densities $\rho_i$ are positive,  the column vector
${\bf v}_i$ is given by $(1,x_i,y_i , \ldots, x_i^n, x_i^{n-1}y_i, \newline \ldots ,$ $x_i y_i^{n-1},y_i^n)^T$,
and the point $(x_i,y_i)$ is in the algebraic variety $\V$ for all $i=1,\ldots,r$.\

In \cite{FiaNie}, L.A. Fialkow and J. Nie  proved abstractly the existence of a representing measure
for the nonsingular quartic MP; they did this using convex analysis techniques. \ They also gave an upper bound of $15$ for the cardinality of the support of a representing measure. \ Using the above mentioned rank-one decomposition, in \cite{CuYoo3} we obtained a  concrete solution of the nonsingular quartic MP, and we were able to discover a method for constructing
a representing (6-atomic) measure using some of the results in \cite{tcmp6}.\
Applying and extending the techniques in \cite{CuYoo3}, we are now able to solve completely the non-extremal sextic
moment problems of rank 7, as follows.

\begin{theorem}\label{main01}
Let $\mathcal V\equiv \mathcal V (\mathcal M(3))$ be the algebraic variety of
$\mathcal M(3)(\beta^{(6)})$ and let $v$ be the cardinality of $\mathcal V$ ($v=+\infty$ is a possible value).\
If $\mathcal M(3)$ is consistent, $\mathcal M(3) \geq 0$, $\mathcal M(2)>0$, $\operatorname{rank \,} \mathcal M(3)=7$, and $v\geq 8$,
then $\beta^{(6)}$ has a $7$-atomic measure.
\end{theorem}

Theorem \ref{main01} covers three cases in Table $1$; it shows that positivity, consistency and the variety condition are often sufficient for solubility.\
However, many instances of TMP as in \cite{CuYoo1},  \cite{CuYoo2}, and \cite{Fia11} show
that a solution to sextic or higher MP requires numerical conditions
associated with the given moment data. \
In general, it is difficult to describe these conditions as generic properties of the moment matrix, so one typically settles for an algorithm. \
This phenomenon is likely the main remaining obstacle in obtaining a general solution
to the other cases listed in Table $1$.\

When specific generic properties of the moment matrix that detect solubility are not available, we find algorithms; for example, in the rank-8 cases with nonzero extremality gaps. \ In our algorithms, we strive to write $\mathcal{M}(n)$ as a sum $\mathcal{M}(n)=\widetilde{ \mathcal{M}(n)} + P$ for some positive moment matrix $P$, which in turn is associated to a few atoms; we then proceed to check if $\widetilde{ \mathcal{M}(n)}$ has a representing measure.

\bigskip

\medskip \textit{Acknowledgments}. \ The authors are deeply grateful to the referee for detecting an error in the original version of the paper, and for making several suggestions that helped improve the presentation. \ Most of the examples, and some of the proofs in this paper, were
obtained using the software tool \textit{Mathematica} \cite{Wol}.


\bigskip


\section{Preliminary Matricial Results}\label{prel}

When we decompose of a moment matrix as a sum $\mathcal{M}(n)=\widetilde{ \mathcal{M}(n)} + P$, the goal is to both reduce the rank (i.e., $\operatorname{rank}\; \widetilde{\mathcal M (n)} < \operatorname{rank}\; \mathcal M (n)$) and obtain a moment matrix $\widetilde{ \mathcal{M}(n)}$ for which we can use previous known results to solve TMP. \ In some cases, we can even make $\widetilde{ \mathcal{M}(n)}$ flat. \  
Controlling the rank of the matrices requires a well known fundamental inequality:

\begin{lemma}\label{rank0} Let $A$ and $B$ be finite matrices. Then
\begin{equation}
\operatorname{rank}\;  (A + B) \leq  \operatorname{rank}\; A +  \operatorname{rank}\; B
\end{equation}
\end{lemma}

Another condition we have to ensure is the positive semidefiniteness of $\widetilde{\mathcal M (n)}$. \ We investigate the eigenvalues of $\widetilde{\mathcal M(n)}$, and need to show its minimal eigenvalue is zero.\ 
We now denote each eigenvalue of an $n\times n$ matrix $A$ as $\lambda_i(A)$ for some $i$ and arrange them
in the ascending order
\begin{equation}
\lambda_1 (A) \leq \cdots\leq \lambda_n(A).
\end{equation}


\subsection{\textbf{Eigenvalue Inequalities}}
The following theorem shows the relationship between the eigenvalues of the matrix and  its perturbation by a rank-one matrix.\

\begin{theorem}\label{rank1}  \cite{HoJo} 
Let $A\in M_{n}$ be Hermitian and let $z\in \mathbb C^n$ be a given vector. If the eigenvalues of $A$ and $A\pm z z^\ast$ are arranged in increasing order as above, we have for $k=1,2,\ldots,n-2$,
\begin{enumerate}[(i)]
\item $\lambda_k(A\pm zz^\ast)\leq \lambda_{k+1}(A)\leq \lambda_{k+2}(A\pm zz^\ast)$,

\item $\lambda_k(A)\leq \lambda_{k+1}(A\pm zz^\ast)\leq  \lambda_{k+2}(A).$
\end{enumerate}
\end{theorem}
Indeed, there is a more general version of the preceding result.\ 

\begin{theorem} \label{rank2} \cite{HoJo}
Let $A,B\in M_n$ be Hermitian and suppose that $B$ has rank at most $r$. Then
\begin{enumerate}[(i)]
\item \label{rank21} $\lambda_k(A+B) \leq \lambda_{k+r}(A) \leq \lambda_{k+2r}(A+B)$ for $k=1,2,\ldots,n-2r$;

\item \label{rank22} $\lambda_k(A) \leq \lambda_{k+r}(A+B) \leq \lambda_{k+2r}(A)$ for $k=1,2,\ldots,n-2r$;

\item \label{rank23}  If $A=U\Lambda U^{\ast}$ with $U=\begin{pmatrix}
{\bf u}_1 \ {\bf u}_2 \ \cdots \ {\bf u}_n\end{pmatrix}\in M_n$ unitary and $\Lambda=\operatorname{diag}(\lambda_1,\ldots,\lambda_n)$ with $\lambda_1\leq \lambda_2\leq \cdots \leq \lambda_n$, and if
$$
B=\lambda_n {\bf u}_n {\bf u}_n^{\ast} + \lambda_{n-1} {\bf u}_{n-1 } {\bf u}_{n-1}^{\ast} +\cdots+ \lambda_{n-r+1} {\bf u}_{n-r+1 } {\bf u}_{n-r+1}^{\ast},
$$
then $\lambda_{\emph{max}}(A-B)=\lambda_{n-r} (A)$.
\end{enumerate}
\end{theorem}

In order to check the positivity of some moment matrices in Section \ref{sec-r8vi}, we will use
the following version of Choleski's Algorithm.\ 

\begin{lemma}
\label{l-asif} \cite{CHO} \ Assume that
\begin{equation*}
P=\left(
\begin{array}{cc}
u & \mathbf{t} \\
\mathbf{t}^{\ast } & P_{0}%
\end{array}%
\right) ,
\end{equation*}%
where $P_{0}$ is an $(n-1)\times (n-1)$ matrix, $\mathbf{t}$ is a row
vector, and $u$ is a real number.\ 
\begin{enumerate}[(i)]
 \item Let $P_{0}$ be invertible. \ Then $\det P=\det P_{0}(u-\mathbf{t}%
P_{0}^{-1}\mathbf{t}^{\ast })$.\ 
\item Let $P_{0}$ be invertible and positive. \ Then 
$$P\geq
0\iff  (u-\mathbf{t}P_{0}^{-1}\mathbf{t}^{\ast })\geq
0\iff \det P\geq 0.$$ 

\item Let $u>0$. \ Then $P\geq 0\iff P_{0}-\mathbf{t}^{\ast
}u^{-1}\mathbf{t}\geq 0$.\ 

\item Let $P\geq 0$ and $p_{ii}=0$ for some $i,$ $1\leq i\leq n$. \ Then $%
p_{ij}=p_{ji}=0$ for all $j=1,\cdots ,n$.
\end{enumerate}
\end{lemma}


\subsection{\textbf{Determinantal Formulas}}

In order to have $\operatorname{rank}\; \widetilde{ \mathcal{M}(n)} < \operatorname{rank}\; \mathcal{M}(n)$ and maintain positivity, we additionally need to study determinantal formulas of the perturbed matrix $\widetilde{ \mathcal{M}(n)}$. \
In the sequel, we will use the following notations: we shall denote by $A=\begin{pmatrix} a_{ij}\end{pmatrix}_{1\leq i,j\leq m}$ an arbitrary square matrix, and by $I_1\equiv I_1(m)$ the $m\times m$ matrix with $1$ in the $(1,1)$-entry and $0$ in all other entries. \  
We shall let $A_{ \{ i_1, \ldots,i_k\} }$ denote the compression of $A$ to the columns and the rows indexed by $\{ i_1, \ldots,i_k\}$. \ 

One easily checks that, for $\alpha\in \mathbb{R}$
\begin{equation}\label{det0}
\det(A-\alpha\, I_1)= \det(A)-\alpha \det \left(A_{\{2,3,\ldots,m\}}\right).
\end{equation}
On the other hand, the Spectral Theorem guarantees that  any rank-one Hermitian matrix is unitarily equivalent to a scalar multiple of $I_1 (m)$. \ If $P$ is rank-one, there exists a unitary operator $U$ such that $U^\ast P U = \lambda\, I_1$, where $\lambda$ is the only nonzero eigenvalue of $P$. \ We now generalize (\ref{det0}) as follows.

\begin{proposition}\label{det1}
Let $A$ be an arbitrary square matrix of size $m$, and let $P$, $U$ and $\lambda$ be as above. \ Then 
\begin{eqnarray*}
\det(A-\rho\, P)= \det(A)-\rho\, \lambda
\det \left((U^\ast AU)_{\{2,3,\cdots,m\}}\right) \; \; \textrm{for all } \rho \in \mathbb{R}.
\end{eqnarray*}
\end{proposition}

\addtocounter{theorem}{-1}

\begin{proof} \ Recall that $U^\ast P U = \lambda\, I_1(m)$. \ Thus,
\begin{eqnarray*}
\det(A-\rho\, P) &=& \det (U^\ast A U -  \rho\, \lambda \, I_1(m)) \\
&=&  \det (U^\ast A U)  -\rho\, \lambda   \det   \left((U^\ast AU)_{\{2,3,\cdots,m\}}\right) \\
&=& \det(A)-\rho\, \lambda
\det \left((U^\ast AU)_{\{2,3,\cdots,m\}}\right).  
\end{eqnarray*}  \qed
\end{proof}


\section{Non-extremal $\mathcal M(3)$ of Rank 7}\label{sec-r7}

The following example of a rank-7, sextic moment problem with an infinite variety was introduced in \cite[Appendix]{Yoo}.\

\begin{example}[An illustration of Theorem \ref{main01}] \label{append}
Beginning with a sextic TMP admitting an $8$-atomic measure, five of whose
atoms lie on the horizontal line $y=1$, we allow one of the moments, $\beta_{60}$, to vary as a parameter.\  The resulting moment matrix $\mathcal M(3)$ is positive semidefinite, it depends on $\beta_{60}$, and it has rank $6$
or $7$ according to the value of $\beta_{60}$: 

\medskip
\noindent\ $\mathcal{M}(3)(\beta_{60}):=$ \newline
\begin{eqnarray*}\small{
\left(
\begin{array}{cccccccccc}
 8 & 19 & 3 & 233 & 35 & 67 & 1441 & 247 & 1 & -303 \\
 19 & 233 & 35 & 1441 & 247 & 1 & 14501 & 1511 & 365 & 455 \\
 3 & 35 & 67 & 247 & 1 & -303 & 1511 & 365 & 455 & 2503 \\
 233 & 1441 & 247 & 14501 & 1511 & 365 & 121489 & 14671 & 1633 & 151 \\
 35 & 247 & 1 & 1511 & 365 & 455 & 14671 & 1633 & 151 & -2111 \\
 67 & 1 & -303 & 365 & 455 & 2503 & 1633 & 151 & -2111 & -16527 \\
 1441 & 14501 & 1511 & 121489 & 14671 & 1633 & \beta_{60} & 122015 & 15245 & 2543 \\
 247 & 1511 & 365 & 14671 & 1633 & 151 & 122015 & 15245 & 2543 & 3413 \\
 1 & 365 & 455 & 1633 & 151 & -2111 & 15245 & 2543 & 3413 & 17615 \\
 -303 & 455 & 2503 & 151 & -2111 & -16527 & 2543 & 3413 & 17615 & 118447
\end{array}
\right) .}
\end{eqnarray*}
\medskip
We are of course interested in the case $\text{rank}\;\mathcal M(3)=7$. \ We can use Gaussian elimination to determine the necessary and sufficient condition for $\text{rank}\; \mathcal M(3)=7$, while
maintaining positive semidefiniteness.\ The moment matrix $\mathcal M(3)(\beta_{60})$ has three column relations (i.e., the columns $X^2Y$, $XY^2$, and $Y^3$ depend on the remaining seven columns), and it is possible to identify precisely the three
polynomials arising from the column relations. \  A calculation shows that the intersection of the corresponding zero sets consists of
the line $y=1$ together with three additional points in the $(x,y)$-plane. \
Thus, the algebraic variety is infinite, while the rank of the moment matrix is $7$. \ $\mathcal M(3)$ is recursively generated and satisfies the variety condition; also, $\mathcal M(3)(\beta_{60})$ admits a flat extension, $\mathcal M(4)$, so that the TMP associated with $\mathcal M(3)$ is soluble, and the $7$-atomic representing measure has support in the algebraic variety. \  Surprisingly, the value of $\beta_{60}$, which played a key role in determining both the rank and the positive semidefiniteness of $\mathcal M(3)$, plays no role in the calculation of the flat extension. \ Thus, the moment matrix automatically generates a family of examples, indexed by $\beta_{60}$. \qed
\end{example}

In contrast with Example \ref{append}, we recall that  there are some moment matrices $\mathcal M(3)$ of
rank 8 or 9, with an infinite variety (see \cite{Fia11} and
\cite{tcmp13}), and which have no representing measure. \ Also, we know that
one can easily find  a pair of cubic polynomials sharing 8 or 9
points, which allow us to construct related moment matrices with
nonzero extremality gaps, with or without representing measures. \ However, it is quite 
difficult to find a triple of cubics intersecting at more than 7
points.\ If one such triple has 7 intersecting points, its associated
moment problem is extremal and a solution was presented in \cite{CuYoo2}, 
using the consistency property. \ 
Nevertheless, we have not been able to find  three cubics with an algebraic variety whose cardinality is 8 or 9.\ Thus, we believe that the sextic MP of rank 7 is more rigid than the corresponding sextic MP with rank 8 or 9. \ Further, one is led to conjecture that positivity, consistency and the variety condition may be sufficient for the existence of a measure, as in Example \ref{append}. \ 

We are now ready to prove Theorem \ref{main01}, which we restate for the reader's convenience.

\begin{theorem}\label{main1}
Let $\mathcal V\equiv \mathcal V (\mathcal M(3))$ be the algebraic variety of
$\mathcal M(3)(\beta^{(6)})$ and let $v$ be the cardinality of $\mathcal V$ ($v=+\infty$ is a possible value).\
If $\mathcal M(3)$ is consistent, $\mathcal M(3) \geq 0$, $\mathcal M(2)>0$, $\operatorname{rank \,} \mathcal M(3)=7$, and $v\geq 8$,
then $\beta^{(6)}$ has a $7$-atomic measure.
\end{theorem}

\addtocounter{theorem}{-1}

\begin{proof} We first claim that there is a point $(a,b)\in\mathcal V$ such that no conic can contain all the points in $\mathcal V - \left\{ (a,b)\right\}$. \ To verify the claim, and without loss of generality, we assume that $v<+\infty$ and we let $(x_1 , y_1), \cdots, (x_v , y_v)$ in $\mathcal V$ denote $v$ distinct points in $\mathcal{V}$. \ We now build the generalized Vandermonde matrix
\begin{eqnarray*}
\mathcal E=\begin{pmatrix}
1 & x_i & y_i & x_i^2 & x_i y_i & y_i^2 & x_i^3 & x_i^2 y_i & x_i y_i^2 & y_i^3
\end{pmatrix}_{i\geq 1}^v.
\end{eqnarray*}
Since the three polynomials associated to the column relations of $\mathcal M(3)$ must pass through all the points in $\mathcal V$, the matrix has only  $7$ linearly independent columns, and a fortiori only $7$ linearly independent rows.\ Thus, we can find a row $R_j$ associated with a point $(x_j,y_j)$ which is a linear combination of the other rows $R_i$ for $i\neq j$; namely, $R_j=\sum_{i\neq j} k_i R_i$ for some $k_i\in \mathbb{R}$.\ Assume now that a quadratic polynomial $c(x,y) \equiv c_1 + c_2 x + c_3 y+ c_4 x^2 + c_5 x y + c_6 y^2$ vanishes on $\mathcal V-\{ (x_j,y_j) \}$, and therefore on every point $(x_i,y_i)$ for $i \ne j$. \ This fact can be described using matrix multiplication, as follows: for $i \ne j$,
$$
R_i \, {\bf a}=\begin{pmatrix}
1& x_i  & y_i& x_i^2  & x_i y_i&  y_i^2 & x_i^3  & x_i^2 y_i& x_i y_i^2 & y_i^3
\end{pmatrix} {\bf a}={\bf 0},
$$
where ${\bf a}=
\begin{pmatrix}
c_1 , c_2  , c_3 , c_4 , c_5  , c_6,0,0,0,0
\end{pmatrix}^T.
$
This leads to
\begin{eqnarray}
R_j \, {\bf a}
= \sum_{i\neq j}k_i R_i \, {\bf a} ={\bf 0}
~ \Longrightarrow ~ c(x_j,y_j)=0.
\end{eqnarray}
Therefore, $c$ also vanishes on $(x_j,y_j)$. \ Since the initial $v$ points were arbitrarily chosen within $\mathcal{V}$, it is clear that $c|_{\mathcal V}\equiv 0$. \ Since $\mathcal{M}(3)$ is consistent, we immediately conclude that $c(X,Y)=\bf{0}$ in the column space of $\mathcal{M}(3)$; but $c$ is a polynomial of degree $2$, so this is a contradiction to the fact that $\mathcal{M}(2)$ is invertible. \ We have thus established the claim.

Now, let us denote by $(a,b)$ the above mentioned point $(x_j,y_j)$, and consider the vector $\textbf v= \left(1,a,b,a^2,ab,b^2,a^3,a^2b,ab^2,b^3
\right)^T$. \ Notice that  ${\textbf v} {\textbf v}^T$ is a rank-one  moment matrix, with the representing measure $\delta_{(a,b)}$. \
We now define
$$
\widetilde {\mathcal M(3)}:= \mathcal M(3) - \rho \, {\textbf v \textbf v}^T,
$$
for some $\rho\neq0$.\
Our strategy is to show that there is a positive $\rho$ such that $\widetilde{ \mathcal M(2)}>0$, $\widetilde {\mathcal M(3)}\geq0$, and $\operatorname{rank \,} \widetilde{ \mathcal M(3)}=6$. \
If so, $\widetilde{ \mathcal M(3)}$ is a flat extension of $\widetilde{ \mathcal M(2)}>0$, and it therefore admits a $6$-atomic representing measure; as a consequence, $\mathcal M(3)$ admits a $7$-atomic representing measure.\ 

To find $\rho$, we first try to bound the rank of $\widetilde{ \mathcal M(3)}$.\
Since $\operatorname{rank}\; \mathcal M(3)\leq \operatorname{rank}\; \widetilde{ \mathcal M(3)} + \operatorname{rank}\; \left({\bf v} {\bf v}^T\right)$ (by Lemma \ref{rank0}), it follows that $\operatorname{rank \,} \widetilde {\mathcal M(3)}\geq 6$.\
On the other hand, the matrix ${\textbf v} {\textbf v}^T$ has the representing measure $\delta_{(a,b)}$ and the point $(a,b)$ is in the  the algebraic variety of $\mathcal M(3)$, and so the three column relations in $\mathcal M(3)$ must hold in ${\textbf v\textbf v}^T$ as well. \
This means $\widetilde{ \mathcal M(3)}$ has at least  three column dependence relations, and therefore $\operatorname{rank \,} \widetilde {\mathcal M(3)}\leq 7$. \ Combining the two estimates, we know that $\operatorname{rank \,} \widetilde{ \mathcal M(3)}$ is 6 or 7.\

We briefly pause to observe that all compressions of the rank-one operator ${\textbf v\textbf v}^T$ to a basis $\mathcal{B}$ of the column space $\mathcal C_{ \mathcal M(3)}$ are still rank-one; indeed the $(1,1)$-entry of ${\textbf v\textbf v}^T$ is $1$, and the column ${\it 1}$ is in $\mathcal{B}$. \ We now return to the search for the desired $\rho$. \ We need to evaluate some determinants, so let $\mathcal{ B}$ be the basis of $\mathcal C_{ \mathcal M(3)}$.\
According to Proposition \ref{det1}, all leading principal minors of $\widetilde {\mathcal M(3)}_{ \mathcal  B}$ are linear in $\rho$. \
In particular,
\begin{eqnarray}\label{eq-r7}
\det\left( \widetilde {\mathcal M(3)}_{ \mathcal  B} \right)=
\det \left(\mathcal M(3)_{ \mathcal  B}\right) - \rho\, \lambda \det \left((U^\ast \mathcal M(3)_{ \mathcal  B} U)_{\{2,3,\ldots,7\}} \right),
\end{eqnarray}
for some unitary matrix $U$, where $\lambda$ is the only nonzero eigenvalue of the rank-one matrix $({\bf v} {\bf v}^T)_{\mathcal{B}}$. \ Since $\mathcal M(3)$ is positive semidefinite, the compression $\mathcal M(3)_{ \mathcal  B}$ is positive definite, which guarantees that both $\det \left(\mathcal M(3)_{ \mathcal  B}\right)$ and $\lambda \det \left((U^\ast \mathcal M(3)_{ \mathcal  B} U)_{\{2,3,\ldots,7\}} \right)$ are positive. \

If we now choose $\rho$ to make the right-hand side of (\ref{eq-r7}) equal to zero, that is, $\rho:= \det \left(\mathcal M(3)_{ \mathcal  B}\right)/ \left(\lambda \det \left((U^\ast \mathcal M(3)_{ \mathcal  B} U)_{\{2,3,\ldots,7\}} \right) \right)$, we see that $\det\left( \widetilde {\mathcal M(3)}_{ \mathcal  B} \right)=0$, and therefore the rank of $\widetilde{ \mathcal M(3)}$ becomes $6$. \
Suppose now that the eigenvalues of $\mathcal M(3)$ and $\widetilde{\mathcal M(3)}$ are arranged in ascending order as in Theorem \ref{rank1}.\
Then we see that $0<\lambda_4(\mathcal M(3))\leq \lambda_5(\widetilde{\mathcal M(3)})$, and hence that $\lambda_k(\widetilde{\mathcal M(3)})>0$ for $k=5,\ldots,10$.\ This tells us that $\widetilde{\mathcal M(3)}$ has exactly 6 positive eigenvalues along with zero whose  multiplicity is 4.\
In other words, $\widetilde{\mathcal M(3)}$ is positive semidefinite.\

Moreover, the proof of the initial claim reveals that the algebraic variety associated with $\widetilde{\mathcal M(3)}$, $\widetilde{\mathcal{V}}$, is obtained from $\mathcal{V}$ by removing the point $(a,b)$. \ This immediately shows that $\widetilde {\mathcal M(2)}$ is positive definite: for, if $\operatorname{rank} \; \widetilde {\mathcal M(2)}<6$, then there would exist a nonzero quadratic form $c(x,y)$ with $c(X,Y)=\bf{0}$ in the column space of $\widetilde {\mathcal M(2)}$; by the Extension Principle \cite[Proposition 3.9]{tcmp1}, we would then conclude that $c(X,Y)=\bf{0}$ in the column space of $\widetilde {\mathcal M(3)}$, and therefore $\widetilde{\mathcal{V}} \subseteq \mathcal{Z}(c)$, that is, the quadratic form $c$ would vanish on $\mathcal{V} - \{(a,b)\}$, a contradiction to the original claim. \ We have thus proved that $\widetilde{\mathcal M(3)}$ is a flat extension of $\widetilde {\mathcal M(2)}$, and it therefore has a $6$-atomic representing measure.\ We conclude that $ \mathcal M(3)$ has a $7$-atomic representing measure, as desired. \qed
\end{proof}



Careful analysis of the proof of Theorem \ref{main1} shows that the following generalization holds.

\begin{corollary}
Let $\mathcal V\equiv \mathcal V (\mathcal{M}(n))$ be the algebraic variety of $\mathcal{M}(n)\left(\beta^{(2n)}\right)$ and let $v$ be the cardinality of $\mathcal V$.\
Suppose $\mathcal{M}(n)$ is consistent, and that $\mathcal{M}(n-1)>0$ and its associated moment sequence has an $r$-atomic representing measure. \ If $\mathcal{M}(n)(\beta^{(2n)})\geq 0$, $\operatorname{rank \,} \mathcal{M}(n)=\frac{n(n+1)}{2}+1$, and $v\geq \frac{n(n+1)}{2}+2$, then $\beta^{(2n)}$  has a (r+1)-atomic measure.\
\end{corollary}


\section{$\mathcal M(3)$ with $r=8$ and $v=9$}\label{sec-r8v9}

To date, most concrete solutions of sextic moment problems include numerical conditions on one or more of the moments; it is generally intricate to express these numerical conditions as specific  properties of the moment matrix.\ Moreover, when we solve a recursively determinate \cite{tcmp13} sextic moment problem ($r=8$ and $v\geq 8$), we need to maintain recursiveness to build the extension $\mathcal M(4)$, and we must also verify the positivity of $\mathcal M(4)$. \ This leads naturally to an algorithmic approach to TMP, as we will see in the following main result. \ We first pause to formulate the key mathematical problem of this section. \

\begin{problem}\label{r8v9} Suppose $\mathcal M(3)\geq 0$ is of rank 8, consistent, with $\mathcal M(2)>0$, and with $v=9$. \
Let $\mathcal V \equiv \{(x_i,y_i)\}_{i=1}^9$ be the algebraic variety of $\mathcal M(3)$. \ Under what conditions does the moment sequence admit a representing measure?
\end{problem}

\bigskip

\begin{algorithm} \label{alg15} \ This Algorithm provides a solution to Problem \ref{r8v9}. 

\medskip
\noindent\textit{Step 1.} \ Build the generalized Vandermonde matrix of $\mathcal V$, namely,
\begin{equation}
\mathcal E:=\begin{pmatrix}
1 & x_i & y_i & x_i^2 & x_i y_i & y_i^2 & x_i^3 & x_i^2 y_i
& x_i y_i^2 & y_i^3
\end{pmatrix}_{i=1}^9.
\end{equation}
Now label the columns of $\mathcal E$ with monomials just as we did for a moment matrix. \ 
Since $\mathcal E$ has 8 linearly independent rows, we can pick a point $(a,b)\in \mathcal V$ such that the row $R_{(a,b)}$ associated with $(a,b)$ is linearly dependent on the other 8 rows.\ 

\medskip
\noindent\textit{Step 2.} \ Let $\mathcal B$ be the basis for the column  space of $\mathcal E$ and let $\mathcal E_{\mathcal B}$ denote the resulting matrix after removing the two dependent columns and the row $R_{(a,b)}$ from $\mathcal E$. \ Observe that $\mathcal E_{\mathcal B}$ is a square matrix of size $8\times 8$; we claim that it is invertible. \ For,  if $\mathcal E_{\mathcal B}$ were singular, there would be another cubic vanishing on $\mathcal V-\{(a,b)\}$; by the Cayley-Bacharach Theorem \cite{EGH}, every cubic passing through any eight of the nine points also passes through the ninth point. \ By Consistency, this new cubic gives rise to a column relation in $\mathcal M(3)$, in addition to the two given column relations. \ This is a contradiction with the hypothesis that the rank of $\mathcal M(3)$ is 8. \

\medskip
\noindent\textit{Step 3.} \ Once we have verified that $\mathcal E_{\mathcal B}$ is invertible, we choose another point $(c,d)\in \mathcal V$ ($(c,d) \ne (a,b)$) and eliminate the row $R_{(c,d)}$ associated with $(c,d)$ from $\mathcal E_{\mathcal B}$; we denote this matrix as $\mathcal E_{\mathcal B}'$.\  Note that $\mathcal E_{\mathcal B}'$ has rank 7 and this fact implies that  there is a new cubic polynomial $r(x,y)$ vanishing on $\widehat{\mathcal V}:=\mathcal V-\{(a,b),(c,d)\}$, besides $p(x,y)$ and $q(x,y)$. The polynomial $r(x,y)$ has a significant role in what follows.\

\medskip
\noindent\textit{Step 4.} We will now use a rank-one decomposition of $\mathcal M(3)$, and will try to understand the structure of the decomposition in case a representing measure exists.\ Suppose $\mathcal M(3)$ has a representing measure.\
Then the variety condition forces a measure to be 8- or 9-atomic.\ 
Let us define a vector-valued function ${\bf v}(x,y):=\left(
1,x,y,x^2,xy ,y^2,x^3,x^2 y ,xy^2,y^3\right)^T$.\
Then, we may write
$$
\mathcal M(3)=\widetilde{\mathcal M(3)}+m_1{\bf v}(a,b) {\bf v}(a,b)^T+m_2{\bf v}(c,d) {\bf v}(c,d)^T,
$$
where $m_1$ and $m_2$ are nonnegative (not simultaneously  zero) for $(a,b),(c,d)\in \mathcal V$.\ 
Notice that  the moment matrices ${\bf v}(a,b) {\bf v}(a,b)^T$ and ${\bf v}(c,d) {\bf v}(c,d)^T$, have the representing measures $\delta_{(a,b)}$ and $\delta_{(c,d)}$, respectively. \ Therefore,  in the presence of a measure, we should be able to find a moment matrix $\widetilde{\mathcal M(3)}$ with a 6- or 7-atomic measure (since $\operatorname{rank}\; \widetilde{\mathcal M(3)}=6$ or 7). \ Denote such representing measure by $\tilde{\mu}$, supported in $\widehat{\mathcal V}$, that is, $\operatorname{supp}\; \tilde\mu \subseteq \widehat{\mathcal V}$.\
Since $\widehat{\mathcal V} \subseteq \mathcal Z(r)$, it follows that $\operatorname{supp}\; \tilde\mu \subseteq \mathcal Z(r)$, and therefore $r(X,Y)={\bf 0}$ in $\widetilde{\mathcal M(3)}$, by Proposition \ref{3-1-tcmp1}. \ 
In short, our goal is to find nonnegative $m_1$ and $m_2$ such that $r(X,Y)={\bf 0}$ in $\widetilde{\mathcal M(3)}$.\ 

\medskip
\noindent\textit{Step 5.} 
In order to find such $m_1$ and $m_2$, we need to solve a linear system of 10 equations with the two unknowns $m_1$ and $m_2$.\ 
If no nonnegative solutions exist, $\mathcal M(3)$ does not have a representing measure.\
In the case when a solution does exist, we must check whether $\widetilde{\mathcal M(3)} \geq 0$ with the fixed $m_1$ and $m_2$ (equivalently, $\Lambda_{\widetilde{\mathcal M(3)}} (x^i y^j r)=0$ for $0\leq i+j\leq 3$).\  

\medskip
\noindent\textit{Step 6.} \ After checking positive semidefiniteness, we still have the two possible cases based on the values of $\operatorname{rank}\;  \widetilde{\mathcal M(3)}$:
If $\operatorname{rank}\; \widetilde{\mathcal M(3)}=6$, then $\widetilde{\mathcal M(3)}$ is a flat extension of $\widetilde{\mathcal M(2)}$; hence $\widetilde{\mathcal M(3)}$ has a 6-atomic  measure, and so $\mathcal M(3)$ has an 8-atomic measure.\ 
Finally, to cover the case  $\operatorname{rank}\; \widetilde{\mathcal M(3)}=7$, notice that $\operatorname{card}\; \mathcal V(\widetilde{\mathcal M(3)}) = 7$; if the cardinality of the variety is 7, then $\widetilde{\mathcal M(3)}$ is extremal and we need to use the results in \cite{CuYoo2}. If $\operatorname{card}\; \mathcal V(\widetilde{\mathcal M(3)})\geq 8$, then it  follows from Theorem \ref{main1} that  $\widetilde{\mathcal M(3)}$ admits a representing measure and so does $\mathcal M(3)$. 

The construction of the Algorithm is therefore complete. \qed
\end{algorithm}

\bigskip

The following example explains how the algorithm works.\

\begin{example}
Consider $\mathcal M(3)$ with two free moments $\beta_{30}$ and $\beta_{40}$:
\begin{eqnarray*}\small{
\left(
\begin{array}{cccccccccc}
 1 & 0 & 0 & 1 & 0 & 1 & \beta_{30} & 0 & -\beta_{30} & 0 \\
 0 & 1 & 0 & \beta_{30} & 0 & -\beta_{30} & \beta_{40} & 0 & \beta_{22} & 0 \\
 0 & 0 & 1 & 0 & -\beta_{30} & 0 & 0 & \beta_{22} & 0 & 9 \\
 1 & \beta_{30} & 0 & \beta_{40} & 0 & \beta_{22} & 41 \beta_{30} & 0 & -16 \beta_{30} & 0 \\
 0 & 0 & -\beta_{30} & 0 & \beta_{22} & 0 & 0 & -16 \beta_{30} & 0 & -9 \beta_{30} \\
 1 & -\beta_{30} & 0 & \beta_{22} & 0 & 9 & -16 \beta_{30} & 0 & -9 \beta_{30} & 0 \\
 \beta_{30} & \beta_{40} & 0 & 41 \beta_{30} & 0 & -16 \beta_{30} & \beta_{60} & 0 & \beta_{42} & 0 \\
 0 & 0 & \beta_{22} & 0 & -16 \beta_{30} & 0 & 0 & \beta_{42} & 0 & \beta_{24} \\
 -\beta_{30} & \beta_{22} & 0 & -16 \beta_{30} & 0 & -9 \beta_{30} & \beta_{42} & 0 & \beta_{24} & 0 \\
 0 & 0 & 9 & 0 & -9 \beta_{30} & 0 & 0 & \beta_{24} & 0 & 81
\end{array}
\right), }
\end{eqnarray*}
where $\beta_{22}=25-\beta_{40}$,  $\beta_{60}=-400+41 \beta_{40}$, $\beta_{42}=400-16 \beta_{40}$, and $\beta_{24}=225-9 \beta_{40}$.
After row reduction, we see that $\mathcal M(3)$ has two column relations:
\begin{eqnarray*}
XY^2=25X -XY \qquad \textrm{and} \qquad Y^3=9Y.\
\end{eqnarray*}
The algebraic variety of these two polynomials is the 9-point set
$$
\mathcal V:=\{(-5,0),(-4,-3),(-4,3),(0,-3),(0,0),(0,3),(4,-3),(4,3),(5,0)\}.
$$
Moreover, using \cite[Proposition 3.6]{Fia08}, a {\it Mathematica} calculation reveals that $\mathcal{M}(3)$ is consistent. \ Let $(a,b):=(4,3)$ and $(c,d):=(5,0)$.\
Then the generalized Vandermonde matrix of $\widehat{\mathcal{V}}:=\mathcal V-\{(a,b),(c,d)\}$ has rank 7.\
Row reduction of the Vandermonde matrix without the two rows associated with the  points $(a,b)$ and $(c,d)$ introduces a new cubic 
$r(x,y):=x^2 y+\frac{1}{3}x^3+4 x y+3 x^2+\frac{20}{3}x$
vanishing on $\mathcal V -\{(a,b),(c,d)\}$.
We write
$$
\mathcal M(3)=\widetilde{\mathcal M(3)}+m_1{\bf v}(a,b) {\bf v}(a,b)^T+m_2{\bf v}(c,d) {\bf v}(c,d)^T
$$
for some $m_1,m_2\in \mathbb R$.\
We now find $m_1$ and $m_2$ such that the column relation $r(X,Y)={\bf 0}$ holds for $\widetilde{\mathcal M(3)}$; this reduces to a system of 10 equation in the two unknowns $m_1$ and $m_2$.\
We can find the solution as $m_1= \frac{1}{576} \left(25-4 \beta_{30}-\beta_{40}\right)$ and $m_2= \frac{1}{450}\left(-16+5 \beta_{30}\right.$ $\left.+\beta_{40}\right)$ for any $\beta_{30}$ and $\beta_{40}$.\
Now, let us fix $\beta_{40}=20$.\
Then a calculation shows that if  $\mathcal M(3)\geq 0$ and $\operatorname{rank}\;\mathcal M(3)=8$, then we need to have $-k <\beta_{30}<k$, where $k:=\sqrt{\frac{3727}{1128}-\frac{\sqrt{7754209}}{1128}}\approx 0.914017$; for positivity of $\widetilde{\mathcal M(3)}$, we should get
$-\frac{5}{4}<\beta_{30}<\frac{4}{5}$;
and for $m_1$ and $m_2$ to be nonzero, it is necessary that
$-\frac{4}{5}\leq\beta_{30}\leq \frac{5}{4}$.\
The common sub-interval for $\beta_{30}$ is then $[-\frac{4}{5},\frac{4}{5}]$. \ We have:
\begin{enumerate}[(i)]
\item If $\beta_{30}=\frac{4}{5}$, then $m_1=\frac{1}{320}$, $m_2=\frac{4}{225}$, $\widetilde{\mathcal M(3)}\geq 0$, and $\operatorname{rank}\; \widetilde{\mathcal M(3)}=6$, which implies  $\mathcal M(3)$ has an 8-atomic measure;

\item If $\beta_{30}=-\frac{4}{5}$, then $m_1=\frac{41}{2880}$, $m_2=0$, $\widetilde{\mathcal M(3)}\geq 0$, and $\operatorname{rank}\; \widetilde{\mathcal M(3)}=7$.\ 
Using \cite[Theorem 6.3.75]{Yoo}, we show  that  $\widetilde{\mathcal M(3)}$ is consistent; this requires finding the quartic polynomial $x^4 +5x^3-16x^2-8x$ vanishing on $\widehat{\mathcal V}$ by looking at the extended, generalized Vandermonde matrix and then testing $\Lambda(x^i y^j (x^4 +5x^3-16x^2-8x))=0$ for $0\leq i+j\leq 2$.\ Since $\widetilde{\mathcal M(3)}$ has a 7-atomic measure, it follows that $\mathcal M(3)$ admits an 8-atomic measure;

\item If $-\frac{4}{5}< \beta_{30}<\frac{4}{5}$, then we know $\mathcal M(3)$ has a 9-atomic measure by checking consistency of $\widetilde{\mathcal M(3)}$ as in (ii).

\item If $-k< \beta_{30}<-\frac{4}{5}$, then both $\mathcal M(3)$ and $\widetilde{\mathcal M(3)}$ are positive but there are no desired $m_1$ and $m_2$;

\item If $\frac{4}{5}< \beta_{30}<k$, then $\mathcal M(3)\geq 0$ but  $\widetilde{\mathcal M(3)}$ is not positive.
\end{enumerate}
Finally, it is possible to show that the last two cases do not admit a representing measure. \qed
\end{example}


\section{$\mathcal M(3)$ with $r=8$ and $v=\infty$}\label{sec-r8vi}

We begin this section by introducing a result that covers singular moment problems with a linear column relation.\
This theorem will be used in Case 2 of Algorithm \ref{Alg}.

\begin{theorem} \label{2-1-tcmp2}
 (cf. \cite[Theorem 2.1]{tcmp2}) Assume that $\mathcal{M}(n)\geq 0$ satisfies (RG) and that $Y=A \textit{1} +B X$ for some $A,B\in \mathbb R$.\
Then $\mathcal{M}(n)$ admits a flat extension $\mathcal{M}(n+1)$.\
\end{theorem}

Theorem \ref{2-1-tcmp2} says that once we have a linear column relation, positivity and recursiveness solve the problem.\
We now consider the last case of $\mathcal M(3)$ with $r=8$.\

\begin{problem}\label{8-inf} Let $\mathcal V$ be the algebraic variety of $\mathcal M(3)$.\   Assume $\mathcal M(3)\geq 0$, of rank 8, consistent,  with $\mathcal M(2)>0$, and with $v=\infty$.\
Under what conditions, does the moment sequence admits a representing measure?
\end{problem}

Since $\operatorname{rank}\; \mathcal M(3)=8$, the moment matrix must have two column dependence relations, say, $p(X,Y)=\bf{0}$ and $q(X,Y)=\bf{0}$.\
In addition, for  $\mathcal V$ to be an infinite set,  $p$ and $q$ must be both reducible and have a common factor.\
If the common factor is a conic, then the  two different line factors of $p$ and $q$ must have only one intersecting point.\
For, if they don't intersect, then all the points in $\mathcal{V}$ are in the zero set of the common  conic and this makes $\mathcal M(2)$  singular.\
On the other hand, if $p$ and $q$ have a common linear factor, then $p$ and $q$ have at least 3 and at most 4 common  points that are not collinear on the graph of the non-common factors of $p$ and $q$.\ (The factors are a pair of different conics or a conic and two lines.)\
If this is not the case, then all the points in the variety can stay in  two lines and this forces to be a  conic column relation in $\mathcal M(3)$.\

Before we describe an algorithm for Problem \ref{8-inf}, we wish to discuss how one finds a solution to  $\mathcal M(3)$ with a conic column relation.\
Due to the equivalence of TMP under the degree-one transformations (see  \cite[Section 5]{tcmp6}), it is sufficient to consider only 5 basic types of nontrivial conics. \
Suppose $\mathcal M(3)$ has rank $7$, is recursively generated, and it has only one \textit{conic} column relation $c(X,Y)=\bf{0}$. \ (Of course, two more column relations must be found in $\mathcal M(3)$, since $\operatorname{rank}\; \mathcal M(3)=7$.)\
Assume also that $\mathcal M(1)>0$.\
First, if the conic $c$ is a parabola or a hyperbola, then $\mathcal M(3)$ admits a representing measure as in \cite{tcmp7} and \cite{tcmp9}.\
Second, if the conic $c$ is an ellipse, then  Theorem 3.5 in \cite{tcmp5} guarantees that $\mathcal M(3)$ has a measure.\
Finally, if the conic $c$ is a pair of intersecting  lines (resp. parallel lines), then we may assume, via a degree-one transformation, that $c(X,Y)=XY$ or $c(X,Y)=X^2-X$. \ The following two propositions show that both cases admit a minimal ($\operatorname{rank}\; \mathcal M(3)$-atomic) representing measure.\ To accomplish this, we will use a separation-of-atoms technique that splits the atoms into two different sets, according to whether they lie in one or the other line. \ These two results will be used in Algorithm \ref{Alg}.

\begin{proposition}\label{xy-0}
Let $\mathcal M(3)$ be a positive semidefinite, recursively generated moment matrix satisfying $XY=0$. \ Then $\mathcal M(3)$ has a 7-atomic representing measure.\
\end{proposition}

\addtocounter{theorem}{-1}

\begin{proof}
Without loss of generality, we can take $\beta_{00}=1$. \ Since $XY=0$, we may write
\begin{equation}
\mathcal M(3)=
\left(
\begin{array}{cccccccccc}
 1 & \beta_{10} & \beta_{01} & \beta_{20} & 0 & \beta_{02} & \beta_{30} & 0 & 0 & \beta_{03} \\
 \beta_{10} & \beta_{20} & 0 & \beta_{30} & 0 & 0 & \beta_{40} & 0 & 0 & 0 \\
 \beta_{01} & 0 & \beta_{02} & 0 & 0 & \beta_{03} & 0 & 0 & 0 & \beta_{04} \\
 \beta_{20} & \beta_{30} & 0 & \beta_{40} & 0 & 0 & \beta_{50} & 0 & 0 & 0 \\
 0 & 0 & 0 & 0 & 0 & 0 & 0 & 0 & 0 & 0 \\
 \beta_{02} & 0 & \beta_{03} & 0 & 0 & \beta_{04} & 0 & 0 & 0 & \beta_{05} \\
 \beta_{30} & \beta_{40} & 0 & \beta_{50} & 0 & 0 & \beta_{60} & 0 & 0 & 0 \\
 0 & 0 & 0 & 0 & 0 & 0 & 0 & 0 & 0 & 0 \\
 0 & 0 & 0 & 0 & 0 & 0 & 0 & 0 & 0 & 0 \\
 \beta_{03} & 0 & \beta_{04} & 0 & 0 & \beta_{05} & 0 & 0 & 0 & \beta_{06}
\end{array}
\right).
\end{equation}
If $\mathcal M(3)$ admits a representing measure, then it admits a finitely atomic representing measure, by the main result in \cite{BaTe}. \ Since the support of such representing measure would be a subset of the degenerate hyperbola $xy=0$, the support can be written as $\{ (x_1,0),\ldots, (x_{\ell_1},0), (0,y_1),\ldots,$ $(0,y_{\ell_2}) \}$, for some nonnegative integers $\ell_1,\ell_2$.\ We now use the above mentioned separation-of-atoms technique, which in this case amounts to decomposing the moment matrix $\mathcal{M}(3)$ as the sum of two moment matrices, each incorporating an additional column relation (and its recursive multiples), either $X=\bf 0$ or $Y=\bf 0$, as follows: 
\begin{eqnarray*}
\mathcal M(3)&=&\left(
\begin{array}{cccccccccc}
 w & 0 & \beta_{01} & 0 & 0 & \beta_{02} & 0 & 0 & 0 & \beta_{03} \\
 0 & 0 & 0 & 0 & 0 & 0 & 0 & 0 & 0 & 0 \\
 \beta_{01} & 0 & \beta_{02} & 0 & 0 & \beta_{03} & 0 & 0 & 0 & \beta_{04} \\
 0 & 0 & 0 & 0 & 0 & 0 & 0 & 0 & 0 & 0 \\
 0 & 0 & 0 & 0 & 0 & 0 & 0 & 0 & 0 & 0 \\
 \beta_{02} & 0 & \beta_{03} & 0 & 0 & \beta_{04} & 0 & 0 & 0 & \beta_{05} \\
 0 & 0 & 0 & 0 & 0 & 0 & 0 & 0 & 0 & 0 \\
 0 & 0 & 0 & 0 & 0 & 0 & 0 & 0 & 0 & 0 \\
 0 & 0 & 0 & 0 & 0 & 0 & 0 & 0 & 0 & 0 \\
 \beta_{03} & 0 & \beta_{04} & 0 & 0 & \beta_{05} & 0 & 0 & 0 & \beta_{06}
\end{array}
\right)
\\
&+&
\left(
\begin{array}{cccccccccc}
 1-w & \beta_{10} & 0 & \beta_{20} & 0 & 0 & \beta_{30} & 0 & 0 & 0 \\
 \beta_{10} & \beta_{20} & 0 & \beta_{30} & 0 & 0 & \beta_{40} & 0 & 0 & 0 \\
 0 & 0 & 0 & 0 & 0 & 0 & 0 & 0 & 0 & 0 \\
 \beta_{20} & \beta_{30} & 0 & \beta_{40} & 0 & 0 & \beta_{50} & 0 & 0 & 0 \\
 0 & 0 & 0 & 0 & 0 & 0 & 0 & 0 & 0 & 0 \\
 0 & 0 & 0 & 0 & 0 & 0 & 0 & 0 & 0 & 0 \\
 \beta_{30} & \beta_{40} & 0 & \beta_{50} & 0 & 0 & \beta_{60} & 0 & 0 & 0 \\
 0 & 0 & 0 & 0 & 0 & 0 & 0 & 0 & 0 & 0 \\
 0 & 0 & 0 & 0 & 0 & 0 & 0 & 0 & 0 & 0 \\
 0 & 0 & 0 & 0 & 0 & 0 & 0 & 0 & 0 & 0
\end{array}
\right),
\end{eqnarray*}
where $w$ is some real number.\
We denote the first matrix in the above sum as $A$ and the second as $B$.\
We attempt to show that both $A$ and $B$  have a representing measure.\
Note first that both moment matrices are recursively generated with a linear column relation.\
All we need to do, by Theorem \ref{2-1-tcmp2}, is to check that the two matrices are positive semidefinite simultaneously for some $w$.\
Consider the rearrangement of the compressions of the two matrices: ${A}_{\{10,6,3,1\}}$ and $ {B}_{\{7,4,2,1\}}$.\
Since the leading principal minors, up to third order, of the two matrices are exactly the same as those of $\mathcal M(3)$, it follows from Lemma \ref{l-asif}(ii) that  
\begin{eqnarray}
A \geq 0, ~ B \geq 0 \iff \det A \geq 0, ~\det B \geq 0
\end{eqnarray}
for some $w.$\
Let $d_A$ (resp. $d_B$) be the third-order leading principal of $A$ (resp. $B$); let $d_6$ be the the sixth-order leading principal of $\mathcal M(3)$.\
Then we can see that 
\begin{eqnarray}\label{eq-alq0}
\det A \geq 0, ~\det B \geq 0
&\iff& \frac{q_A}{d_A}\leq w  \leq   \frac{q_B}{d_B}
\end{eqnarray}
where if $w=q_A/d_A$ (respectively $w=q_B/d_B$), then $\det A = 0$ (respectively, $\det B = 0$).\
Here we have a fortuitous coincidence; a calculation shows that
\begin{eqnarray}
\frac{q_B}{d_B} - \frac{q_A}{d_A} = \frac{d_6}{d_A d_B} > 0,
\end{eqnarray}
from which we conclude that  there is always a desired $w$ satisfying (\ref{eq-alq0}).\ In particular, if we take $w$ as one of the two end points $q_A/d_A$ or $q_B/d_A$, then one of $A$ and $B$ has rank 3 (in other words, $A$ or $B$ is flat).\  Therefore, $\mathcal M(3)$ admits a 7-atomic measure. \qed
\end{proof}

In similar fashion, we can cover the remaining case:

\begin{proposition}\label{x2-x}
Let $\mathcal M(3)$ be a positive semidefinite, recursively generated moment matrix, with column relation $X^2=X$. \ Then $\mathcal{M}(3)$ admits a $7$-atomic representing measure.\
\end{proposition}

\addtocounter{theorem}{-1}

\begin{proof}
The condition $X^2=X$ allows us to write
\begin{equation}
\mathcal M(3)=\left(
\begin{array}{cccccccccc}
 \beta_{00} & \beta_{10} & \beta_{01} & \beta_{10} & \beta_{11} & \beta_{02} & \beta_{10} & \beta_{11} & \beta_{12} & \beta_{03} \\
 \beta_{10} & \beta_{10} & \beta_{11} & \beta_{10} & \beta_{11} & \beta_{12} & \beta_{10} & \beta_{11} & \beta_{12} & \beta_{13} \\
 \beta_{01} & \beta_{11} & \beta_{02} & \beta_{11} & \beta_{12} & \beta_{03} & \beta_{11} & \beta_{12} & \beta_{13} & \beta_{04} \\
 \beta_{10} & \beta_{10} & \beta_{11} & \beta_{10} & \beta_{11} & \beta_{12} & \beta_{10} & \beta_{11} & \beta_{12} & \beta_{13} \\
 \beta_{11} & \beta_{11} & \beta_{12} & \beta_{11} & \beta_{12} & \beta_{13} & \beta_{11} & \beta_{12} & \beta_{13} & \beta_{14} \\
 \beta_{02} & \beta_{12} & \beta_{03} & \beta_{12} & \beta_{13} & \beta_{04} & \beta_{12} & \beta_{13} & \beta_{14} & \beta_{05} \\
 \beta_{10} & \beta_{10} & \beta_{11} & \beta_{10} & \beta_{11} & \beta_{12} & \beta_{10} & \beta_{11} & \beta_{12} & \beta_{13} \\
 \beta_{11} & \beta_{11} & \beta_{12} & \beta_{11} & \beta_{12} & \beta_{13} & \beta_{11} & \beta_{12} & \beta_{13} & \beta_{14} \\
 \beta_{12} & \beta_{12} & \beta_{13} & \beta_{12} & \beta_{13} & \beta_{14} & \beta_{12} & \beta_{13} & \beta_{14} & \beta_{15} \\
 \beta_{03} & \beta_{13} & \beta_{04} & \beta_{13} & \beta_{14} & \beta_{05} & \beta_{13} & \beta_{14} & \beta_{15} & \beta_{06}
\end{array}
\right).
\end{equation}
In the presence of a measure for  $\mathcal M(3)$,
we should be able to separate moments (except $\beta_{00}$) into two groups by the location of atoms; in detail, if  $\{ (x_1,0),\ldots, (x_{\ell_1},0),$ $(x_{\ell_1 +1},1), (x_{\ell_2 +2},1),\ldots,(x_{\ell_1+\ell_2},1) \}$ is the support of a measure for some integers $\ell_1,\ell_2$, then we can write
\begin{eqnarray}\label{eq-x2}
\mathcal M(3)=A_0 + A_1\equiv
\begin{pmatrix} \alpha_{ij}^{(0)}\end{pmatrix}+
\begin{pmatrix} \alpha_{ij}^{(1)}\end{pmatrix},
\end{eqnarray}
where $A_0$ (respectively, $A_1$) is the moment matrix formed by the atoms on the line $x=0$ (respectively, $x=1$).\
We next observe that the moments $\beta_{ij}$ ($0<j\leq6$) are generated by the atoms lying only on $x=1$ so that  $\beta_{ij}=\alpha_{ij}^{(1)}$ ($0<j\leq6$).\
Since $A_1$ is recursively generated, it must have the column relations: $X^2=1$, $XY=Y$, $X^3=1$, $X^2Y=Y$, and $X Y^2=Y^2$.\
These additional relations determine the rest of moments in $A_1$ except $\alpha_{06}^{(1)}$; indeed, we can readily show that $\alpha_{1j}^{(1)}=\beta_{1j}$ for $0\leq j\leq 5$.\
We now let $\alpha_{06}^{(1)}:=\tau$ (look for $\tau$ in the lower right-hand corner of each matrix summand in the next displayed matrix equation). \ It follows that the sum in (\ref{eq-x2}) becomes 
\begin{eqnarray}
\mathcal M(3)&=&\left(
\begin{array}{cccccccccc}
 \beta_{00}-\beta_{10} & 0 & \beta_{01}-\beta_{11} & 0 & 0 & \beta_{02}-\beta_{12} & 0 & 0 & 0 & \beta_{03}-\beta_{13} \\
 0 & 0 & 0 & 0 & 0 & 0 & 0 & 0 & 0 & 0 \\
 \beta_{01}-\beta_{11} & 0 & \beta_{02}-\beta_{12} & 0 & 0 & \beta_{03}-\beta_{13} & 0 & 0 & 0 & \beta_{04}-\beta_{14} \\
 0 & 0 & 0 & 0 & 0 & 0 & 0 & 0 & 0 & 0 \\
 0 & 0 & 0 & 0 & 0 & 0 & 0 & 0 & 0 & 0 \\
 \beta_{02}-\beta_{12} & 0 & \beta_{03}-\beta_{13} & 0 & 0 & \beta_{04}-\beta_{14} & 0 & 0 & 0 & \beta_{05}-\beta_{15} \\
 0 & 0 & 0 & 0 & 0 & 0 & 0 & 0 & 0 & 0 \\
 0 & 0 & 0 & 0 & 0 & 0 & 0 & 0 & 0 & 0 \\
 0 & 0 & 0 & 0 & 0 & 0 & 0 & 0 & 0 & 0 \\
 \beta_{03}-\beta_{13} & 0 & \beta_{04}-\beta_{14} & 0 & 0 & \beta_{05}-\beta_{15} & 0 & 0 & 0 & \beta_{06}-\tau
\end{array}
\right)  \nonumber \\
&+&
\left(
\begin{array}{cccccccccc}
 \beta_{10} & \beta_{10} & \beta_{11} & \beta_{10} & \beta_{11} & \beta_{12} & \beta_{10} & \beta_{11} & \beta_{12} & \beta_{13} \\
 \beta_{10} & \beta_{10} & \beta_{11} & \beta_{10} & \beta_{11} & \beta_{12} & \beta_{10} & \beta_{11} & \beta_{12} & \beta_{13} \\
 \beta_{11} & \beta_{11} & \beta_{12} & \beta_{11} & \beta_{12} & \beta_{13} & \beta_{11} & \beta_{12} & \beta_{13} & \beta_{14} \\
 \beta_{10} & \beta_{10} & \beta_{11} & \beta_{10} & \beta_{11} & \beta_{12} & \beta_{10} & \beta_{11} & \beta_{12} & \beta_{13} \\
 \beta_{11} & \beta_{11} & \beta_{12} & \beta_{11} & \beta_{12} & \beta_{13} & \beta_{11} & \beta_{12} & \beta_{13} & \beta_{14} \\
 \beta_{12} & \beta_{12} & \beta_{13} & \beta_{12} & \beta_{13} & \beta_{14} & \beta_{12} & \beta_{13} & \beta_{14} & \beta_{15} \\
 \beta_{10} & \beta_{10} & \beta_{11} & \beta_{10} & \beta_{11} & \beta_{12} & \beta_{10} & \beta_{11} & \beta_{12} & \beta_{13} \\
 \beta_{11} & \beta_{11} & \beta_{12} & \beta_{11} & \beta_{12} & \beta_{13} & \beta_{11} & \beta_{12} & \beta_{13} & \beta_{14} \\
 \beta_{12} & \beta_{12} & \beta_{13} & \beta_{12} & \beta_{13} & \beta_{14} & \beta_{12} & \beta_{13} & \beta_{14} & \beta_{15} \\
 \beta_{13} & \beta_{13} & \beta_{14} & \beta_{13} & \beta_{14} & \beta_{15} & \beta_{13} & \beta_{14} & \beta_{15} & \tau
\end{array}
\right).
\end{eqnarray}
We now observe that $\mathcal M(3)$ has a representing measure if and only if so do $A_0$ and $A_1$, and this is equivalent to $A_0\geq0$ and $A_1\geq 0$ for some $\alpha>0$ (by Theorem \ref{2-1-tcmp2}).\
To check the positivity, let $d_J$ be the leading principal minors of the compression of $\mathcal M(3)$ to the rows and columns in the index set $J$; also, let $d_k^{(i)}$'s be the $k$-th leading principal minors of the compression of $A_i$ to the basis of its column space for $i=1,2$ and $1\leq k\leq 4$.\
We first notice that all diagonal entries in $A_1$ are positive (except $\tau$) and we can easily check that $d_1^{(1)}= d_{\{2\}}>0$, $d_2^{(1)}= d_{\{2, 5\}}>0$, and  $d_3^{(1)}= d_{\{2, 5,9\}}>0$.\
We next observe that the positivity of the leading principal minors implies that the diagonal entries in $A_0$ are positive (except $\beta_{06}-\tau$).\
Also, a calculation shows that
\begin{eqnarray}
d_2^{(0)}= d_{\{1,2, 3,5\}}/d_2^{(1)}>0,\qquad
d_3^{(0)}= d_{\{1,2, 3,5,6,9\}}/d_3^{(1)}>0.
\end{eqnarray}
Finally, for the positivity of both $A_0$ and $A_1$ we must have 
\begin{eqnarray}
q_1\leq \tau \leq q_0,
\end{eqnarray}
where $q_0$ (respectively, $q_1$) is  a quantity depending on moments such that $d_4^{(0)}=0$ (respectively, $d_4^{(1)}=0$).\
As before, we have a fortuitous coincidence that 
\begin{eqnarray}
q_0 - q_1 = \frac{\det \mathcal M(3)_{\{1,2,3,5,6,8,9\}}}{d_3^{(0)}d_3^{(1)} } > 0.
\end{eqnarray}
To complete the proof, we apply Lemma \ref{l-asif}.  \qed
\end{proof}


\bigskip

We are ready to consider the cases of $r=8$ with an infinite variety.

\begin{algorithm} \label{Alg} \ In this algorithm we provide a solution to Problem \ref{8-inf}. \ Let us write $p(x,y)=\ell_1(x,y) c_1(x,y)$ and $q(x,y)=\ell_2(x,y) c_2(x,y)$, where $\ell_i$ is a line and $c_i$ is a conic for $i=1,2$.\

\medskip
\noindent\textit{Case 1.} $c(x,y) \equiv c_1(x,y)=c_2(x,y)$

Let $ (a_0,b_0)\in \mathcal Z(\ell_1)\cap \mathcal Z (\ell_2)$.\
Then $(a_0,b_0)$ must be in the support of a representing measure; otherwise, consistency of $\mathcal M(3)$ forces to be a conic column relation in $\mathcal M(3)$.\
Notice that $\mathcal M(3)$ has a representing measure if and only if we may write $ \mathcal M(3)$ as a sum of two moment matrices, that is, for some $\rho_0 >0$,
\begin{equation}
\mathcal M(3)=  \rho_0 {\bf v} {\bf v}^T+ \mathcal M_c(3),
\end{equation}
where $\mathbf v =\left(
1,a_0 , b_0, a_0^2,a_0b_0 , b_0^2,a_0^3,a_0^2 b_0 ,a_0b_0^2,b_0^3
\right)^T$
and  $\mathcal M_c(3)$ is a moment matrix generated with atoms in the graph of $c(x,y)=0$.\
Observe that by Proposition \ref{3-1-tcmp1}, $\mathcal M_c (3)$ must be endowed with the quadratic column relation $c(X,Y)=\bf{0}$, and hence $\mathcal M_c (3)$ has at least 3 column relations.\
Indeed, we can claim that there is
 a positive $\rho_0$ such that $\operatorname{rank}\; \mathcal M_c(3)=7$.\
For,  let $\mathcal{ B}$ be the basis of the column space of ${ \mathcal M(3)}$.\ By Proposition \ref{det1}, all leading principal minors of $\mathcal M_c(3)_{ \mathcal  B}$ are linear  in $\rho$; \; that is, 
\begin{eqnarray}
\det\left(  \mathcal M_c(3)_{ \mathcal  B} \right)=
\det \left(\mathcal M(3)_{ \mathcal  B}\right) - \rho\, \lambda \det
\left((U^\ast \mathcal M(3)_{ \mathcal  B} U)_{\{2,3,\ldots,8\}} \right),
\end{eqnarray}
for some unitary matrix $U$, where $\lambda$ is the only nonzero eigenvalue of $({\bf v} {\bf v}^T)_{\mathcal{B}}$.\
Positive definiteness of  $\mathcal M(3)_{ \mathcal  B}$ implies that both $\det \left(\mathcal M(3)_{ \mathcal  B}\right)$ and $\lambda \det \left((U^\ast \mathcal M(3)_{ \mathcal  B} U)_{\{2,3,\ldots,8\}} \right)$ are positive.\
We thus can take $\rho_0$ as $\det \left(\mathcal M(3)_{ \mathcal  B}\right)/ \left(\lambda \det \left((U^\ast \mathcal M(3)_{ \mathcal  B} U)_{\{2,3,\ldots,8\}} \right) \right)$.\  
Next, if the eigenvalues of $\mathcal M(3)$ and $\mathcal M_c(3)$ are arranged in ascending order as in Theorem \ref{rank1}, then we can see that $0<\lambda_3(\mathcal M(3))\leq \lambda_4(\widetilde {\mathcal M_c(3)}) $.\
Since $\operatorname{rank}\; \mathcal M_c(3)= 7$, it follows that $\mathcal M_c(3)$ has the eigenvalue zero with multiplicity 3, through which we can conclude that
$ \lambda_k(\mathcal M_c(3))=0$ for $k=1,2,3$ and
$ \lambda_k(\mathcal M_c(3))>0$ for $k=4,\ldots,10$.\ In other words, $\mathcal M_c(3)$ is positive semidefinite.\

In summary, with the specific $\rho_0$, $\mathcal M(3)$ has a representing measure if and only if $\mathcal M_c(3)$ has a representing measure.\
If $c$ is a circle, parabola, or hyperbola, then we know $\mathcal M_c(3)$ admits a representing measure; if $c$ is a pair of intersecting or parallel lines, then we need to  apply Propositions \ref{xy-0} or \ref{x2-x}.\

\medskip
\noindent\textit{Case 2.} $\ell(x,y) \equiv \ell_1(x,y)=\ell_2(x,y)$

Let $ (c_i,d_i)\in \mathcal Z(c_1)\cap \mathcal Z (c_2) $ for $i=1,\ldots, m_2$ ($3\leq m_2\leq 4$).\ 
Similarly, $\mathcal M(3)$ has a representing measure if and only if  $ \mathcal M(3)$ can be written as a sum of two moment matrices:
\begin{equation}
\mathcal M(3)={\mathcal M_\ell (3)}+
\begin{pmatrix} \beta^{(c)}_{ij} \end{pmatrix}
\equiv {\mathcal M_\ell (3)}+\mathcal M_c(3),
\end{equation}
where 
$ \beta^{(c)} _{ij} = \sum_{k=1}^{m_2} \rho_k^{(c)} c_k^i d_k^j $ for some positive $\rho_i^{(c)}$ ($1\leq i \leq m_2$) and ${\mathcal M_\ell (3)}$ is a moment matrix generated by atoms in the line $\ell$.\ 
We next need to see  that ${\mathcal M_\ell (3)}$ must have the column relation $\ell(X,Y)=\bf{0}$.\ 
Applying the relation to $\mathcal M_\ell(3)=\mathcal M(3)-\mathcal M_c(3)$, we have a linear system of $10$ equations in the unknowns, $ \rho_1^{(c)},\ldots, \rho_{m_2}^{(c)} $ (at least 3 of them are positive).\ 
If the system does not have a nonnegative solution set, then $\mathcal M(3)$ does not have a representing measure.\
If we can find a solution of the system,   then since $\mathcal M_c(3)$ obviously has a representing measure, it follows from Theorem \ref{2-1-tcmp2} that  we just need to check if $\mathcal M_\ell(3) \geq 0$ and if $\mathcal M_\ell(3)$ satisfies (RG).\ 
If $\mathcal M_\ell(3)$ passes both tests, then it has a representing measure and consequently, so does $\mathcal M(3)$. \qed
\end{algorithm}






%


\end{document}